# Certain Subclass of Harmonic Multivalent Functions Defined by New Linear Operator

**Ali H. Maran [1], Abdul Rahman S. Juma [2], Raheam A. Al-Saphory [1]***

[1]Department of Mathematics, College of Education for Pure Sciences, University of Tikrit, Tikrit, IRAQ.

[2]Department of Mathematics, College of Education for Pure Sciences, University of Anbar, Ramadi, IRAQ.

*Corresponding Author: Raheam A. Al-Saphory

**ABSTRACT:** The main goal of the present paper is to introduce a new class of harmonic multivalent functions defined by a new linear operator in the open unit disc $\mathbb{U} = \{z: |z| < 1\}$. Thus, some geometric properties have examined, including coefficient inequality, extreme points, convolution conditions, convex linear combinations, and integral transforms for the class $\mathcal{H}_{p,q,\ell}^{\delta}(p, q, \Upsilon, \ell, \delta)$.

## 1. INTRODUCTION

The study of analytic and harmonic functions is of great importance in the field of mathematics and its applications [1]. These functions are essential in complex analysis, for example the theory of geometric functions, and harmonic mapping, they play a crucial role in various mathematical contexts, such as the study of Riemann surfaces and the theory of univalent functions, which make them a powerful concept in the field of mathematics [2]. One of their most important properties is starlikeness, convexity, close to convexity, extreme points, and the convolution condition [3][4].

A continuous function $f = u + iv$ is complex valued harmonic function in a complex domain $\mathbb{C}$, where both $u$ and $v$ real harmonic in $\mathbb{C}$ [5]. Thus, in any simply connected domain $D \subseteq \mathbb{C}$ [6],[7].We can write $f = h + \bar{g}$ where $h$ and $g$ are analytic in $D$ [8],[9].We call $h$ the analytic part and $g$ the co-analytic part of $f$, a necessary and sufficient condition for $f$ to be locally univalent and sense-preserving in $D$ is that $|h'(z)| > |g'(z)|$ in $D$ [10]. Denoted by $\mathcal{H}_{p,q,m}^{\delta}(p, q, \Upsilon, \ell, \delta)$ the class of function $f = h + \bar{g}$ that are harmonic multivalent and sense pressing in the unit disk $\mathbb{U} = \{z: |z| < 1\}$, for $f = h + \bar{g} \in \mathcal{H}_{p,q,\ell}^{\delta}(p, q, \Upsilon, \ell, \delta)$, we may express the analytic functions $h$ and $g$ as [11],[12] by the form

$$h(z) = z^{\ell} + \sum_{\Upsilon=\ell+1}^{\infty} a_{\Upsilon} z^{\Upsilon}, \quad g(z) = \sum_{\Upsilon=\ell}^{\infty} b_{\Upsilon} z^{\Upsilon}, |b_{\ell}| < 1. \tag{1}$$

Now, Sharma *et al.* [13] study the following operator

$$\vartheta_{\ell}(p, q, \delta; z) = z^{\ell} + \sum_{\Upsilon=\ell+1}^{\infty} \frac{\left([\delta+\ell]_{p,q}\right)_{\Upsilon-\ell}}{[\Upsilon-\ell]_{p,q}} z^{\Upsilon}, \tag{2}$$

where, $\delta > -\ell, \ell \in \mathbb{N}, 0 < q < p \leq 1$.
and Hussain *et al.* [14] investigate the following operator

$$\varphi(\ell, q, t, z) = z^{\ell} + \sum_{\Upsilon=\ell+1}^{\infty} [\Upsilon + (\ell - 1)]_q^t z^{\Upsilon}, \tag{3}$$

Here, the Hadamard product (or convolution) of the above mappings $\vartheta_{\ell}(p, q, \delta; z)$ and $\varphi(\ell, q, t, z)$ is characterized the new linear operator via the next form

$$\mathcal{H}_{p,q,\ell}^{\delta} = (\vartheta_{\ell} * \varphi) f(z) = z^{\ell} + \sum_{\Upsilon=\ell+1}^{\infty} \frac{[\Upsilon+(\ell-1)]_q^t \left([\delta+\ell]_{p,q}\right)_{\Upsilon-\ell}}{[\Upsilon-\ell]_{p,q}} a_{\Upsilon} z^{\Upsilon}. \tag{4}$$

*Corresponding author: saphory@tu.edu.iq

Thus,

$$\mathcal{H}^{\delta}_{p,q,\ell}f(z)=\mathcal{H}^{\delta}_{p,q,\ell}h(z)+\overline{\mathcal{H}^{\delta}_{p,q,\ell}g(z)},$$

where

$$\mathcal{H}^{\delta}_{p,q,\ell}h(z)=z^{\ell}+\sum_{\Upsilon=\ell+1}^{\infty}\frac{[\Upsilon+(\ell-1)]^{t}_{q}\left([\delta+\ell]_{p,q}\right)_{\Upsilon-\ell}}{[\Upsilon-\ell]_{p,q}}a_{\Upsilon}z^{\Upsilon}.$$

and

$$\overline{\mathcal{H}^{\delta}_{p,q,\ell}g(z)}=\sum_{\Upsilon=\ell}^{\infty}\frac{[\Upsilon+(\ell-1)]^{t}_{q}\left([\delta+\ell]_{p,q}\right)_{\Upsilon-\ell}}{[\Upsilon-\ell]_{p,q}}\overline{b_{\Upsilon}z^{\Upsilon}}.$$

Hence

$$\mathcal{H}^{\delta}_{p,q,\ell}f(z)=z^{\ell}+\sum_{\Upsilon=\ell+1}^{\infty}\frac{[\Upsilon+(\ell-1)]^{t}_{q}\left([\delta+\ell]_{p,q}\right)_{\Upsilon-\ell}}{[\Upsilon-\ell]_{p,q}}a_{\Upsilon}z^{\Upsilon}+\sum_{\Upsilon=\ell}^{\infty}\frac{[\Upsilon+(\ell-1)]^{t}_{q}\left([\delta+\ell]_{p,q}\right)_{\Upsilon-\ell}}{[\Upsilon-\ell]_{p,q}}\overline{b_{\Upsilon}z^{\Upsilon}}. \tag{5}$$

Thus, one can define a subclass of $\mathcal{H}^{\delta}_{p,q,\ell}(p,q,\Upsilon,\ell,\delta)$ multivalent harmonic functions, which satisfying the conditions:

$$R\left\{\frac{z^{2}\left(\mathcal{H}^{\delta}_{p,q,\ell}h(z)\right)''-\overline{z^{2}\left(\mathcal{H}^{\delta}_{p,q,\ell}g(z)\right)''}}{z\left(\mathcal{H}^{\delta}_{p,q,\ell}h(z)\right)'+\overline{z\left(\mathcal{H}^{\delta}_{p,q,\ell}g(z)\right)'}}\right\}\geq\sigma, \tag{6}$$

where $0\leq\sigma<1,\delta>-\ell,\ell\in\mathbb{N},0<q<p\leq 1$.

**Lemma 1.1.[15]** Let $\alpha\geq 0$. Then $R\{\theta\}\geq\alpha$ if and only if $|\theta-(1+\alpha)|\geq|\theta+(1-\alpha)|$, where $\theta$ be any complex number.

## 2. COEFFICIENT INEQUALITY

In this section, we prove that the sufficient condition for $f$ to be in the class $\mathcal{H}^{\delta}_{p,q,\ell}(p,q,\Upsilon,\ell,\delta)$, given by the next outcome.

**Theorem 2.1.** Let $f=h+\overline{g}$ is define by (1). If

$$[\Upsilon(\Upsilon-1)+\Upsilon(1-\sigma)]\sum_{\Upsilon=\ell+1}^{\infty}\frac{[\Upsilon+(\ell-1)]^{t}_{q}\left([\delta+\ell]_{p,q}\right)_{\Upsilon-\ell}}{[\Upsilon-\ell]_{p,q}}a_{\Upsilon}z^{\Upsilon}+[\Upsilon(\Upsilon-1)-\Upsilon(1+\sigma)]\sum_{\Upsilon=\ell}^{\infty}\frac{[\Upsilon+(\ell-1)]^{t}_{q}\left([\delta+\ell]_{p,q}\right)_{\Upsilon-\ell}}{[\Upsilon-\ell]_{p,q}}\overline{b_{\Upsilon}z^{\Upsilon}}\geq 0. \tag{7}$$

Then, $f$ is harmonic multivalent sense – preserving $\mathbb{U}$ and $\in\mathcal{H}^{\delta}_{p,q,\ell}(p,q,\Upsilon,\ell,\delta)$ .

**Proof**. To prove that if the inequality (7) holds, then

$$R\left\{\frac{z^{2}\left(\mathcal{H}^{\delta}_{p,q,\ell}h(z)\right)''-\overline{z^{2}\left(\mathcal{H}^{\delta}_{p,q,\ell}g(z)\right)''}}{z\left(\mathcal{H}^{\delta}_{p,q,\ell}h(z)\right)'+\overline{z\left(\mathcal{H}^{\delta}_{p,q,\ell}g(z)\right)'}}\right\}=Re\left\{\frac{A(z)}{B(z)}\right\},$$

where

$$A(z)=z^{2}\left(\mathcal{H}^{\delta}_{p,q,\ell}h(z)\right)''-\overline{z^{2}\left(\mathcal{H}^{\delta}_{p,q,\ell}g(z)\right)''},$$

and

$$B(z)=z\left(\mathcal{H}^{\delta}_{p,q,\ell}h(z)\right)'+\overline{z\left(\mathcal{H}^{\delta}_{p,q,\ell}g(z)\right)'}$$

By using lemma 1.1 it can get

$$|A(z)+(1-\sigma)B(z)|-|A(z)-(1+\sigma)B(z)|\geq 0 \tag{8}$$

Then, by the substituted $A(z)$ and $B(z)$ in the inequality (1.8), we obtain

$$\left|z^{2}\left(\mathcal{H}^{\delta}_{p,q,\ell}h(z)\right)''-\overline{z^{2}\left(\mathcal{H}^{\delta}_{p,q,\ell}g(z)\right)''}+(1-\sigma)\left(z\left(\mathcal{H}^{\delta}_{p,q,\ell}h(z)\right)'+\overline{z\left(\mathcal{H}^{\delta}_{p,q,\ell}g(z)\right)'}\right)\right|-\left|z^{2}\left(\mathcal{H}^{\delta}_{p,q,\ell}h(z)\right)''-\overline{z^{2}\left(\mathcal{H}^{\delta}_{p,q,\ell}g(z)\right)''}-(1+\sigma)\left(z\left(\mathcal{H}^{\delta}_{p,q,\ell}h(z)\right)'+\overline{z\left(\mathcal{H}^{\delta}_{p,q,\ell}g(z)\right)'}\right)\right|\geq 0$$

Hence,

$$\left|\ell(\ell-1)z^{\ell}+\Upsilon(\Upsilon-1)\sum_{\Upsilon=\ell+1}^{\infty}\frac{[\Upsilon+(\ell-1)]^{t}_{q}\left([\delta+\ell]_{p,q}\right)_{\Upsilon-\ell}}{[\Upsilon-\ell]_{p,q}}a_{\Upsilon}z^{\Upsilon}-\Upsilon(\Upsilon-1)\sum_{\Upsilon=\ell}^{\infty}\frac{[\Upsilon+(\ell-1)]^{t}_{q}\left([\delta+\ell]_{p,q}\right)_{\Upsilon-\ell}}{[\Upsilon-\ell]_{p,q}}\overline{b_{\Upsilon}z^{\Upsilon}}+(1-\sigma)\left(\ell z^{\ell}+\Upsilon\sum_{\Upsilon=\ell+1}^{\infty}\frac{[\Upsilon+(\ell-1)]^{t}_{q}\left([\delta+\ell]_{p,q}\right)_{\Upsilon-\ell}}{[\Upsilon-\ell]_{p,q}}a_{\Upsilon}z^{\Upsilon}+\Upsilon\sum_{\Upsilon=\ell}^{\infty}\frac{[\Upsilon+(\ell-1)]^{t}_{q}\left([\delta+\ell]_{p,q}\right)_{\Upsilon-\ell}}{[\Upsilon-\ell]_{p,q}}\overline{b_{\Upsilon}z^{\Upsilon}}\right)\right|-\left|\ell(\ell-1)z^{\ell}+\Upsilon(\Upsilon-1)\sum_{\Upsilon=\ell+1}^{\infty}\frac{[\Upsilon+(\ell-1)]^{t}_{q}\left([\delta+\ell]_{p,q}\right)_{\Upsilon-\ell}}{[\Upsilon-\ell]_{p,q}}a_{\Upsilon}z^{\Upsilon}-\Upsilon(\Upsilon-1)\sum_{\Upsilon=\ell}^{\infty}\frac{[\Upsilon+(\ell-1)]^{t}_{q}\left([\delta+\ell]_{p,q}\right)_{\Upsilon-\ell}}{[\Upsilon-\ell]_{p,q}}\overline{b_{\Upsilon}z^{\Upsilon}}-(1+\sigma)\left(\ell z^{\ell}+\Upsilon\sum_{\Upsilon=\ell+1}^{\infty}\frac{[\Upsilon+(\ell-1)]^{t}_{q}\left([\delta+\ell]_{p,q}\right)_{\Upsilon-\ell}}{[\Upsilon-\ell]_{p,q}}a_{\Upsilon}z^{\Upsilon}+\Upsilon\sum_{\Upsilon=\ell}^{\infty}\frac{[\Upsilon+(\ell-1)]^{t}_{q}\left([\delta+\ell]_{p,q}\right)_{\Upsilon-\ell}}{[\Upsilon-\ell]_{p,q}}\overline{b_{\Upsilon}z^{\Upsilon}}\right)\right|\geq 0.$$

Consequently,

$\Big|[\ell(\ell-1)+\ell-\ell\sigma]z^{\ell}+\curlyvee(\curlyvee-1)\sum_{\curlyvee=\ell+1}^{\infty}\frac{[\curlyvee+(\ell-1)]_q^t([\delta+\ell]_{p,q})_{\curlyvee-\ell}}{[\curlyvee-\ell]_{p,q}}a_{\curlyvee}z^{\curlyvee}-\curlyvee(\curlyvee-1)\sum_{\curlyvee=\ell}^{\infty}\frac{[\curlyvee+(\ell-1)]_q^t([\delta+\ell]_{p,q})_{\curlyvee-\ell}}{[\curlyvee-\ell]_{p,q}}\overline{b_{\curlyvee}z^{\curlyvee}}+$ $\ell z^{\ell}+\curlyvee\sum_{\curlyvee=\ell+1}^{\infty}\frac{[\curlyvee+(\ell-1)]_q^t([\delta+\ell]_{p,q})_{\curlyvee-\ell}}{[\curlyvee-\ell]_{p,q}}a_{\curlyvee}z^{\curlyvee}+\curlyvee\sum_{\curlyvee=\ell}^{\infty}\frac{[\curlyvee+(\ell-1)]_q^t([\delta+\ell]_{p,q})_{\curlyvee-\ell}}{[\curlyvee-\ell]_{p,q}}\overline{b_{\curlyvee}z^{\curlyvee}}-\sigma\ell z^{\ell}-$ $\sigma\curlyvee\sum_{\curlyvee=\ell+1}^{\infty}\frac{[\curlyvee+(\ell-1)]_q^t([\delta+\ell]_{p,q})_{\curlyvee-\ell}}{[\curlyvee-\ell]_{p,q}}a_{\curlyvee}z^{\curlyvee}-\sigma\curlyvee\sum_{\curlyvee=\ell}^{\infty}\frac{[\curlyvee+(\ell-1)]_q^t([\delta+\ell]_{p,q})_{\curlyvee-\ell}}{[\curlyvee-\ell]_{p,q}}\overline{b_{\curlyvee}z^{\curlyvee}}\Big|-\Big|\ell[(\ell-1)-\ell-\ell\sigma]z^{\ell}+\curlyvee(\curlyvee-$ $1)\sum_{\curlyvee=\ell+1}^{\infty}\frac{[\curlyvee+(\ell-1)]_q^t([\delta+\ell]_{p,q})_{\curlyvee-\ell}}{[\curlyvee-\ell]_{p,q}}a_{\curlyvee}z^{\curlyvee}-\curlyvee(\curlyvee-1)\sum_{\curlyvee=\ell}^{\infty}\frac{[\curlyvee+(\ell-1)]_q^t([\delta+\ell]_{p,q})_{\curlyvee-\ell}}{[\curlyvee-\ell]_{p,q}}\overline{b_{\curlyvee}z^{\curlyvee}}+\ell z^{\ell}-$ $\curlyvee\sum_{\curlyvee=\ell+1}^{\infty}\frac{[\curlyvee+(\ell-1)]_q^t([\delta+\ell]_{p,q})_{\curlyvee-\ell}}{[\curlyvee-\ell]_{p,q}}a_{\curlyvee}z^{\curlyvee}-\curlyvee\sum_{\curlyvee=\ell}^{\infty}\frac{[\curlyvee+(\ell-1)]_q^t([\delta+\ell]_{p,q})_{\curlyvee-\ell}}{[\curlyvee-\ell]_{p,q}}\overline{b_{\curlyvee}z^{\curlyvee}}-\sigma\ell z^{\ell}-$ $\sigma\curlyvee\sum_{\curlyvee=\ell+1}^{\infty}\frac{[\curlyvee+(\ell-1)]_q^t([\delta+\ell]_{p,q})_{\curlyvee-\ell}}{[\curlyvee-\ell]_{p,q}}a_{\curlyvee}z^{\curlyvee}-\sigma\curlyvee\sum_{\curlyvee=\ell}^{\infty}\frac{[\curlyvee+(\ell-1)]_q^t([\delta+\ell]_{p,q})_{\curlyvee-\ell}}{[\curlyvee-\ell]_{p,q}}\overline{b_{\curlyvee}z^{\curlyvee}}\Big|$

$$\geq 2[\ell(\ell-2-\sigma)]|z|^{\curlyvee}-2[\curlyvee(\curlyvee-1)+\curlyvee(1-\sigma)]\sum_{\curlyvee=\ell+1}^{\infty}\frac{[\curlyvee+(\ell-1)]_q^t([\delta+\ell]_{p,q})_{\curlyvee-\ell}}{[\curlyvee-\ell]_{p,q}}|a_{\curlyvee}||z|^{\curlyvee}$$
$$-2[\curlyvee(\curlyvee-1)-\curlyvee(1+\sigma)]\sum_{\curlyvee=\ell}^{\infty}\frac{[\curlyvee+(\ell-1)]_q^t([\delta+\ell]_{p,q})_{\curlyvee-\ell}}{[\curlyvee-\ell]_{p,q}}\overline{|b_{\curlyvee}||z|^{\curlyvee}}\geq 0.$$

By virtue of (7), this implies that $f\in\mathcal{H}_{p,q,\ell}^{\delta}(p,q,\curlyvee,\ell,\delta)$. ∎

Thus, the following theorem demonstrates the condition (6) and it is also necessary for functions $f(z)=h(z)+\overline{g(z)}$.

**Theorem 2.2.** Suppose $f\in\mathcal{H}_{p,q,\ell}^{\delta}(p,q,\curlyvee,\ell,\delta)$, if and only if

$$[\curlyvee(\curlyvee-1)+\curlyvee(1-\sigma)]\sum_{\curlyvee=\ell+1}^{\infty}\frac{[\curlyvee+(\ell-1)]_q^t([\delta+\ell]_{p,q})_{\curlyvee-\ell}}{[\curlyvee-\ell]_{p,q}}a_{\curlyvee}z^{\curlyvee}+[\curlyvee(\curlyvee-1)-\curlyvee(1+\sigma)]\sum_{\curlyvee=\ell}^{\infty}\frac{[\curlyvee+(\ell-1)]_q^t([\delta+\ell]_{p,q})_{\curlyvee-\ell}}{[\curlyvee-\ell]_{p,q}}\overline{b_{\curlyvee}z^{\curlyvee}}\geq 0. \quad (9)$$

**Proof.** we only need to prove the " only if " of the theorem. to this end, for functions $f$ of the form (1) with condition (6), we notice that the condition

$$R\left\{\frac{z^2\left(\mathcal{H}_{p,q,\ell}^{\delta}h(z)\right)''-\overline{z^2\left(\mathcal{H}_{p,q,\ell}^{\delta}g(z)\right)''}}{z\left(\mathcal{H}_{p,q,\ell}^{\delta}h(z)\right)'+\overline{z\left(\mathcal{H}_{p,q,\ell}^{\delta}g(z)\right)'}}\right\}\geq\sigma,$$

The above inequality is equivalent to

$$R\left\{\frac{(\ell-1)z^{\ell}+\curlyvee(\curlyvee-1)\sum_{\curlyvee=\ell+1}^{\infty}\frac{[\curlyvee+(\ell-1)]_q^t([\delta+\ell]_{p,q})_{\curlyvee-\ell}}{[\curlyvee-\ell]_{p,q}}a_{\curlyvee}z^{\curlyvee}-\curlyvee(\curlyvee-1)\sum_{\curlyvee=\ell}^{\infty}\frac{[\curlyvee+(\ell-1)]_q^t([\delta+\ell]_{p,q})_{\curlyvee-\ell}}{[\curlyvee-\ell]_{p,q}}\overline{b_{\curlyvee}z^{\curlyvee}}}{\ell z^{\ell}+\curlyvee\sum_{\curlyvee=\ell+1}^{\infty}\frac{[\curlyvee+(\ell-1)]_q^t([\delta+\ell]_{p,q})_{\curlyvee-\ell}}{[\curlyvee-\ell]_{p,q}}a_{\curlyvee}z^{\curlyvee}+\curlyvee\sum_{\curlyvee=\ell}^{\infty}\frac{[\curlyvee+(\ell-1)]_q^t([\delta+\ell]_{p,q})_{\curlyvee-\ell}}{[\curlyvee-\ell]_{p,q}}\overline{b_{\curlyvee}z^{\curlyvee}}}\right\}\geq 0 \quad (10)$$

Thus, this inequality must hold all values of $z$, such that $|z|=r<1$.

So, the above inequality can reduce to

$$\frac{(\ell-1)+\curlyvee(\curlyvee-1)\sum_{\curlyvee=\ell+1}^{\infty}\frac{[\curlyvee+(\ell-1)]_q^t([\delta+\ell]_{p,q})_{\curlyvee-\ell}}{[\curlyvee-\ell]_{p,q}}a_{\curlyvee}-\curlyvee(\curlyvee-1)\sum_{\curlyvee=\ell}^{\infty}\frac{[\curlyvee+(\ell-1)]_q^t([\delta+\ell]_{p,q})_{\curlyvee-\ell}}{[\curlyvee-\ell]_{p,q}}b_{\curlyvee}}{\ell+\curlyvee\sum_{\curlyvee=\ell+1}^{\infty}\frac{[\curlyvee+(\ell-1)]_q^t([\delta+\ell]_{p,q})_{\curlyvee-\ell}}{[\curlyvee-\ell]_{p,q}}a_{\curlyvee}+\curlyvee\sum_{\curlyvee=\ell}^{\infty}\frac{[\curlyvee+(\ell-1)]_q^t([\delta+\ell]_{p,q})_{\curlyvee-\ell}}{[\curlyvee-\ell]_{p,q}}b_{\curlyvee}}\geq 0 \quad (11)$$

In the case where $r\to 1^{-}$ and if the condition in (9) doesn't hold, then the numerator in (11) is negative. This is a contradiction with our assumptions and then, $f\in\mathcal{H}_{p,q,\ell}^{\delta}(p,q,\curlyvee,\ell,\delta)$.

Since $f\in\mathcal{H}_{p,q,\ell}^{\delta}(p,q,\curlyvee,\ell,\delta)$, then for

$$\sum_{\curlyvee=\ell+1}^{\infty}|X_{\curlyvee}|+\sum_{\curlyvee=\ell}^{\infty}|Y_{\curlyvee}|=1$$

The harmonic univalent function

$$f(z)=z^{\ell}+\sum_{\curlyvee=\ell+1}^{\infty}\frac{\ell(\ell-2-\sigma)[\curlyvee-\ell]_{p,q}}{[\curlyvee(\curlyvee-1)+\curlyvee(1-\sigma)][\curlyvee+(\ell-1)]_q^t([\delta+\ell]_{p,q})_{\curlyvee-\ell}}X_{\curlyvee}z^{\curlyvee}$$
$$+\sum_{\curlyvee=\ell}^{\infty}\frac{\ell(\ell-2-\sigma)[\curlyvee-\ell]_{p,q}}{[\curlyvee(\curlyvee-1)-\curlyvee(1+\sigma)][\curlyvee+(\ell-1)]_q^t([\delta+\ell]_{p,q})_{\curlyvee-\ell}}\overline{Y_{\curlyvee}z^{\curlyvee}}.$$

Demonstrates the equality in the coefficient bound given by (9) is sharp.

## 3. EXTREME POINT

Now, the next result proves extreme point in the class $\mathcal{H}_{p,q,\ell}^{\delta}(p,q,\curlyvee,\ell,\delta)$.

**Theorem 3.1.** A function $f=h+\overline{g}\in\mathcal{H}_{p,q,\ell}^{\delta}(p,q,\curlyvee,\ell,\delta)$ if and only if $f(z)$ can be expressed in the form

$$f(z)=\sum_{\Upsilon=\ell}^{\infty}[X_\Upsilon h(z)+Y_\Upsilon g(z)], \tag{12}$$

where

$$h_\ell(z)=z^\ell$$

and

$$h(z)=z^\ell+\frac{\ell(\ell-2-\sigma)[\Upsilon-\ell]_{p,q}}{[\Upsilon(\Upsilon-1)+\Upsilon(1-\sigma)][\Upsilon+(\ell-1)]_q^t\left([\delta+\ell]_{p,q}\right)_{\Upsilon-\ell}}z^\Upsilon, \Upsilon\geq 2.$$

with,

$$g(z)=z^\ell+\frac{\ell(\ell-2-\sigma)[\Upsilon-\ell]_{p,q}}{[\Upsilon(\Upsilon-1)-\Upsilon(1+\sigma)][\Upsilon+(\ell-1)]_q^t\left([\delta+\ell]_{p,q}\right)_{\Upsilon-\ell}}\overline{z^\Upsilon}, \Upsilon\geq 1$$

So, $h(z)\geq 0$ , $g(z)\geq 0$ .

$$X_\ell=1-\sum_{\Upsilon=\ell+1}^{\infty}X_\Upsilon+\sum_{\Upsilon=\ell}^{\infty}Y_\Upsilon$$

The extreme points of $f\in\mathcal{H}_{p,q,\ell}^{\delta}(p,q,\Upsilon,\ell,\delta)$ are $h(z)$ and $g(z)$ .

**Proof.** Let $f(z)$ is define by (12). Then, we have

$$f(z)=\sum_{\Upsilon=\ell}^{\infty}[X_\Upsilon h(z)+Y_\Upsilon g(z)]$$

Hence,

$$f(z)=\sum_{\Upsilon=\ell}^{\infty}(X_\Upsilon+Y_\Upsilon)z^\ell+\sum_{\Upsilon=\ell+1}^{\infty}\frac{\ell(\ell-2-\sigma)[\Upsilon-\ell]_{p,q}}{[\Upsilon(\Upsilon-1)+\Upsilon(1-\sigma)][\Upsilon+(\ell-1)]_q^t\left([\delta+\ell]_{p,q}\right)_{\Upsilon-\ell}}X_\Upsilon z^\ell$$
$$+\sum_{\Upsilon=\ell}^{\infty}\frac{\ell(\ell-2-\sigma)[\Upsilon-\ell]_{p,q}}{[\Upsilon(\Upsilon-1)-\Upsilon(1+\sigma)][\Upsilon+(\ell-1)]_q^t\left([\delta+\ell]_{p,q}\right)_{\Upsilon-\ell}}Y_\Upsilon z^\ell.$$

Thus,

$$f(z)=z^\ell+\sum_{\Upsilon=\ell+1}^{\infty}\frac{\ell(\ell-2-\sigma)[\Upsilon-\ell]_{p,q}}{[\Upsilon(\Upsilon-1)+\Upsilon(1-\sigma)][\Upsilon+(\ell-1)]_q^t\left([\delta+\ell]_{p,q}\right)_{\Upsilon-\ell}}X_\Upsilon z^\ell$$
$$+\sum_{\Upsilon=\ell}^{\infty}\frac{\ell(\ell-2-\sigma)[\Upsilon-\ell]_{p,q}}{[\Upsilon(\Upsilon-1)-\Upsilon(1+\sigma)][\Upsilon+(\ell-1)]_q^t\left([\delta+\ell]_{p,q}\right)_{\Upsilon-\ell}}Y_\Upsilon z^\ell.$$

Furthermore, let

$$f(z)=\sum_{\Upsilon=\ell+1}^{\infty}\frac{[\Upsilon(\Upsilon-1)+\Upsilon(1-\sigma)][\Upsilon+(\ell-1)]_q^t\left([\delta+\ell]_{p,q}\right)_{\Upsilon-\ell}}{\ell(\ell-2-\sigma)[\Upsilon-\ell]_{p,q}}a_\Upsilon z^\ell$$
$$+\sum_{\Upsilon=\ell}^{\infty}\frac{[\Upsilon(\Upsilon-1)-\Upsilon(1+\sigma)][\Upsilon+(\ell-1)]_q^t\left([\delta+\ell]_{p,q}\right)_{\Upsilon-\ell}}{\ell(\ell-2-\sigma)[\Upsilon-\ell]_{p,q}}b_\Upsilon z^\ell.$$

Then,

$$f(z)=\sum_{\Upsilon=\ell+1}^{\infty}\frac{[\Upsilon(\Upsilon-1)+\Upsilon(1-\sigma)][\Upsilon+(\ell-1)]_q^t\left([\delta+\ell]_{p,q}\right)_{\Upsilon-\ell}}{\ell(\ell-2-\sigma)[\Upsilon-\ell]_{p,q}}\times$$
$$\left(\frac{\ell(\ell-2-\sigma)[\Upsilon-\ell]_{p,q}}{[\Upsilon(\Upsilon-1)+\Upsilon(1-\sigma)][\Upsilon+(\ell-1)]_q^t\left([\delta+\ell]_{p,q}\right)_{\Upsilon-\ell}}X_\Upsilon\right)$$
$$+\sum_{\Upsilon=\ell}^{\infty}\frac{[\Upsilon(\Upsilon-1)-\Upsilon(1+\sigma)][\Upsilon+(\ell-1)]_q^t\left([\delta+\ell]_{p,q}\right)_{\Upsilon-\ell}}{\ell(\ell-2-\sigma)[\Upsilon-\ell]_{p,q}}\times$$
$$\left(\frac{\ell(\ell-2-\sigma)[\Upsilon-\ell]_{p,q}}{[\Upsilon(\Upsilon-1)-\Upsilon(1+\sigma)][\Upsilon+(\ell-1)]_q^t\left([\delta+\ell]_{p,q}\right)_{\Upsilon-\ell}}Y_\Upsilon\right)$$

Consequently,

$$f(z)=\sum_{\Upsilon=\ell+1}^{\infty}X_\Upsilon+\sum_{\Upsilon=\ell}^{\infty}Y_\Upsilon=1-X_\ell\leq 1$$

Hence, $f(z)\in\mathcal{H}_{p,q,\ell}^{\delta}(p,q,\Upsilon,\ell,\delta)$.

Conversely, if $f(z)\in\mathcal{H}_{p,q,\ell}^{\delta}(p,q,\Upsilon,\ell,\delta)$ . Assume that for $(\Upsilon\geq 2)$

$$X_\ell=1-\sum_{\Upsilon=\ell+1}^{\infty}X_\Upsilon-\sum_{\Upsilon=\ell}^{\infty}Y_\Upsilon,$$

where

$$X_\Upsilon = \frac{[\Upsilon(\Upsilon-1)+\Upsilon(1-\sigma)][\Upsilon+(\ell-1)]_q^t\left([\delta+\ell]_{p,q}\right)_{\Upsilon-\ell}}{\ell(\ell-2-\sigma)[\Upsilon-\ell]_{p,q}}|a_\Upsilon|.$$

$$Y_\Upsilon = \frac{[\Upsilon(\Upsilon-1)-\Upsilon(1+\sigma)][\Upsilon+(\ell-1)]_q^t\left([\delta+\ell]_{p,q}\right)_{\Upsilon-\ell}}{\ell(\ell-2-\sigma)[\Upsilon-\ell]_{p,q}}|b_\Upsilon|.$$

One can put

$$f(z) = z^\ell + \sum_{\Upsilon=\ell+1}^{\infty} a_\Upsilon z^\Upsilon + \sum_{\Upsilon=\ell}^{\infty} b_\Upsilon \overline{z^\Upsilon}$$

Then,

$$f(z) = z^\ell + \sum_{\Upsilon=\ell+1}^{\infty} \frac{\ell(\ell-2-\sigma)[\Upsilon-\ell]_{p,q}}{[\Upsilon(\Upsilon-1)+\Upsilon(1-\sigma)][\Upsilon+(\ell-1)]_q^t\left([\delta+\ell]_{p,q}\right)_{\Upsilon-\ell}} X_\Upsilon z^\Upsilon + \sum_{\Upsilon=\ell}^{\infty} \frac{\ell(\ell-2-\sigma)[\Upsilon-\ell]_{p,q}}{[\Upsilon(\Upsilon-1)-\Upsilon(1+\sigma)][\Upsilon+(\ell-1)]_q^t\left([\delta+\ell]_{p,q}\right)_{\Upsilon-\ell}} Y_\Upsilon \overline{z^\Upsilon}.$$

Therefore,

$$f(z) = z^\ell + \sum_{\Upsilon=\ell+1}^{\infty} [z^\ell + h(z)]X_\Upsilon + \sum_{\Upsilon=\ell}^{\infty} [z^\ell + g(z)]Y_\Upsilon.$$

Consequently,

$$f(z) = \left[1 - \sum_{\Upsilon=\ell+1}^{\infty} X_\Upsilon - \sum_{\Upsilon=\ell}^{\infty} Y_\Upsilon\right] z^\ell + \sum_{\Upsilon=\ell+1}^{\infty} X_\Upsilon h(z) + \sum_{\Upsilon=\ell}^{\infty} Y_\Upsilon g(z).$$

Then,

$$f(z) = \sum_{\Upsilon=\ell+1}^{\infty} X_\Upsilon h(z) + \sum_{\Upsilon=\ell}^{\infty} Y_\Upsilon g(z). \blacksquare$$

## 4. CONVEX LINEAR COMBINATIONS

Next result demonstrates convex linear combinations for the class $\mathcal{H}_{p,q,\ell}^{\delta}(p,q,\Upsilon,\ell,\delta)$.

**Theorem 4.1.** The class $\mathcal{H}_{p,q,\ell}^{\delta}(p,q,\Upsilon,\ell,\delta)$ is closed convex linear combination.

**Proof.** Let each of the functions

$$f_j(z) = z^\ell + \sum_{\Upsilon=\ell+1}^{\infty} |a_{\Upsilon,j}| z^\Upsilon + \sum_{\Upsilon=\ell}^{\infty} |b_{\Upsilon,j}| \overline{z^\Upsilon}, (a_{\Upsilon,j}, b_{\Upsilon,j} \geq 0; j = 1,2)$$

be in the class $\mathcal{H}_{p,q,\ell}^{\delta}(p,q,\Upsilon,\ell,\delta)$. Thus, it is sufficient to show that the function $G(z)$ defined by

$$G(z) = \mu f_1(z) + (1-\mu) f_2(z), (0 \leq \mu < 1)$$

is also in class $\mathcal{H}_{p,q,\ell}^{\delta}(p,q,\Upsilon,\ell,\delta)$.

Since $(0 \leq \mu \leq 1)$, then

$$G(z) = z^\ell + \sum_{\Upsilon=\ell+1}^{\infty} \left[\mu|a_{\Upsilon,1}| + (1-\mu)|a_{\Upsilon,2}|\right] z^\Upsilon + \sum_{\Upsilon=\ell}^{\infty} \left[\mu|b_{\Upsilon,1}| + (1-\mu)|b_{\Upsilon,2}|\right]\overline{z^\Upsilon},$$

and with the aid of theorem 2.2, implies

$$\sum_{\Upsilon=\ell+1}^{\infty} \left[\mu|a_{\Upsilon,1}| + (1-\mu)|a_{\Upsilon,2}|\right] z^\Upsilon + \sum_{\Upsilon=\ell}^{\infty} \left[\mu|b_{\Upsilon,1}| + (1-\mu)|b_{\Upsilon,2}|\right]\overline{z^\Upsilon}.$$

$$\sum_{\Upsilon=\ell+1}^{\infty} [\Upsilon(\Upsilon-1)+\Upsilon(1-\sigma)]\left[\mu|a_{\Upsilon,1}| + (1-\mu)|a_{\Upsilon,2}|\right] z^\Upsilon + \sum_{\Upsilon=\ell}^{\infty} [\Upsilon(\Upsilon-1)+\Upsilon(1-\sigma)]\left[\mu|b_{\Upsilon,1}| + (1-\mu)|b_{\Upsilon,2}|\right]\overline{z^\Upsilon}.$$

Therefore,

$$= \mu \sum_{\Upsilon=\ell+1}^{\infty} \left([\Upsilon(\Upsilon-1)+\Upsilon(1-\sigma)]|a_{\Upsilon,1}| + [\Upsilon(\Upsilon-1)-\Upsilon(1+\sigma)]|b_{\Upsilon,1}|\right) z^\Upsilon$$

$$+(1-\mu) \sum_{\Upsilon=\ell+1}^{\infty} \left([\Upsilon(\Upsilon-1)+\Upsilon(1-\sigma)]|a_{\Upsilon,2}| + [\Upsilon(\Upsilon-1)-\Upsilon(1+\sigma)]|b_{\Upsilon,2}|\right) z^\Upsilon.$$

$$\leq \mu\ell[\{(\ell-2-\sigma)-(\ell-2+\sigma)\}|b_\ell|] + (1-\mu)\ell[\{(\ell-2-\sigma)-(\ell-2+\sigma)\}|b_\ell|]$$

$$\leq \mu\ell[\{(\ell-2-\sigma)-(\ell-2+\sigma)\}|b_\ell|]$$

which show that $G(z) \in \mathcal{H}^{\delta}_{p,q,\ell}(p,q,\curlyvee,\ell,\delta)$. ■

## 5. CONVOLUTION CONDITION

For harmonic functions

$$f(z) = z^{\ell} + \sum_{\curlyvee=\ell+1}^{\infty}|a_{\curlyvee}|z^{\curlyvee} + \sum_{\curlyvee=\ell}^{\infty}|b_{\curlyvee}|\,\overline{z}^{\curlyvee} \tag{13}$$

and

$$M(z) = z^{\ell} + \sum_{\curlyvee=\ell+1}^{\infty}|c_{\curlyvee}|\,z^{\curlyvee} + \sum_{\curlyvee=\ell}^{\infty}|d_{\curlyvee}|\,\overline{z}^{\curlyvee} \tag{14}$$

Defines the convolution of $f$ and $M$ as

$$(f * M)(z) = z^{\ell} + \sum_{\curlyvee=\ell+1}^{\infty}|a_{\curlyvee}c_{\curlyvee}|\,z^{\curlyvee} + \sum_{\curlyvee=\ell}^{\infty}|b_{\curlyvee}d_{\curlyvee}|\,\overline{z}^{\curlyvee}. \tag{15}$$

The following theorem examines the convolution property of the class $\mathcal{H}^{\delta}_{p,q,\ell}(p,q,\curlyvee,\ell,\delta)$.

**Theorem 5.1.** If $f$ and $M \in \mathcal{H}^{\delta}_{p,q,\ell}(p,q,\curlyvee,\ell,\delta)$, then $f * g \in \mathcal{H}^{\delta}_{p,q,\ell}(p,q,\curlyvee,\ell,\delta)$.

**Proof.** Since $f(z) = z^{\ell} + \sum_{\curlyvee=\ell+1}^{\infty} a_{\curlyvee}\, z^{\curlyvee} + \sum_{\curlyvee=\ell}^{\infty} b_{\curlyvee}\, z^{\curlyvee}$,
be in the class $\mathcal{H}^{\delta}_{p,q,\ell}(p,q,\ell,\delta)$, and

$$g(z) = z^{\ell} + \sum_{\curlyvee=\ell+1}^{\infty} c_{\curlyvee}\, z^{\curlyvee} + \sum_{\curlyvee=\ell}^{\infty} d_{\curlyvee}\, z^{\curlyvee}.$$

be in $\mathcal{H}^{\delta}_{p,q,\ell}(p,q,\curlyvee,\ell,\delta)$.

Consider convolution functions $f * g$ the following:

$$[\curlyvee(\curlyvee-1)+\curlyvee(1-\sigma)]\sum_{\curlyvee=\ell+1}^{\infty}\frac{[\curlyvee+(\ell-1)]_q^t\big([\delta+\ell]_{p,q}\big)_{\curlyvee-\ell}}{[\curlyvee-\ell]_{p,q}}a_{\curlyvee}c_{\curlyvee}z^{\curlyvee}$$
$$+[\curlyvee(\curlyvee-1)\;-\curlyvee(1+\sigma)]\sum_{\curlyvee=\ell}^{\infty}\frac{[\curlyvee+(\ell-1)]_q^t\big([\delta+\ell]_{p,q}\big)_{\curlyvee-\ell}}{[\curlyvee-\ell]_{p,q}}b_{\curlyvee}d_{\curlyvee}z^{\curlyvee}$$
$$\leq[\curlyvee(\curlyvee-1)+\curlyvee(1-\sigma)]\sum_{\curlyvee=\ell+1}^{\infty}\frac{[\curlyvee+(\ell-1)]_q^t\big([\delta+\ell]_{p,q}\big)_{\curlyvee-\ell}}{[\curlyvee-\ell]_{p,q}}a_{\curlyvee}z^{\curlyvee}$$
$$+[\curlyvee(\curlyvee-1)\;-\curlyvee(1+\sigma)]\sum_{\curlyvee=\ell}^{\infty}\frac{[\curlyvee+(\ell-1)]_q^t\big([\delta+\ell]_{p,q}\big)_{\curlyvee-\ell}}{[\curlyvee-\ell]_{p,q}}b_{\curlyvee}z^{\curlyvee}\leq 1. ■$$

## 6. INTEGRAL TRANSFORMS

This section, presents an important result for the class $\mathcal{H}^{\delta}_{p,q,\ell}(p,q,\curlyvee,\ell,\delta)$ under the generalizes Bernardi-Libera-Livingston integral operator which described by

$$F_u\big(f(z)\big) = \frac{u+p}{z^u}\int_0^z t^{u-1}f(z)dt\,, u > -1 \tag{16}$$

**Theorem 6.1.** Let $f_n(z) \in \mathcal{H}^{\delta}_{p,q,\ell}(p,q,\curlyvee,\ell,\delta)$, then $F_u\big(f_{\curlyvee}(z)\big) \in \overline{\mathcal{H}^{\delta}_{p,q,\ell}}(p,q,\curlyvee,\ell,\delta)$ .

**Proof.** From definition of $F_u\big(f_{\curlyvee}(z)\big)$ given by (16), it follows that

$$F_u\big(f_{\curlyvee}(z)\big) = \frac{u+p}{z^u}\int_0^z t^{u-1}\left[t^p + \sum_{\curlyvee=\ell+1}^{\infty} a_{\curlyvee}\, t^{\curlyvee} + \overline{\sum_{\curlyvee=\ell}^{\infty} b_{\curlyvee}\, t^{\curlyvee}}\right]dt.$$
$$= z^p + \frac{u+p}{\curlyvee+u}\sum_{\curlyvee=\ell+1}^{\infty} a_{\curlyvee}\, z^{\curlyvee} + \frac{u+p}{\curlyvee+u}\sum_{\curlyvee=\ell}^{\infty} b_{\curlyvee}\,\overline{z^{\curlyvee}}.$$
$$= z^p + \sum_{\curlyvee=\ell+1}^{\infty}\Psi_{\curlyvee}\, z^{\curlyvee} + \sum_{\curlyvee=\ell}^{\infty}\mathrm{K}_{\curlyvee}\,\overline{z^{\curlyvee}},$$

where,

$$\Psi_{\curlyvee} = \frac{u+p}{\curlyvee+u}a_{\curlyvee} \text{ and } \mathrm{K}_{\curlyvee} = \frac{u+p}{\curlyvee+u}b_{\curlyvee}.$$

Hence,

$$[\curlyvee(\curlyvee-1)+\curlyvee(1-\sigma)]\sum_{\curlyvee=\ell+1}^{\infty}\frac{[\curlyvee+(\ell-1)]_q^t\left([\delta+\ell]_{p,q}\right)_{\curlyvee-\ell}}{[\curlyvee-\ell]_{p,q}}\Psi_\curlyvee$$
$$+[\curlyvee(\curlyvee-1)\;-\curlyvee(1+\sigma)]\sum_{\curlyvee=\ell}^{\infty}\frac{[\curlyvee+(\ell-1)]_q^t\left([\delta+\ell]_{p,q}\right)_{\curlyvee-\ell}}{[\curlyvee-\ell]_{p,q}}\mathrm{K}_\curlyvee$$
$$\leq[\curlyvee(\curlyvee-1)+\curlyvee(1-\sigma)]\sum_{\curlyvee=\ell+1}^{\infty}\frac{[\curlyvee+(\ell-1)]_q^t\left([\delta+\ell]_{p,q}\right)_{\curlyvee-\ell}}{[\curlyvee-\ell]_{p,q}}a_\curlyvee$$
$$+[\curlyvee(\curlyvee-1)\;-\curlyvee(1+\sigma)]\sum_{\curlyvee=\ell}^{\infty}\frac{[\curlyvee+(\ell-1)]_q^t\left([\delta+\ell]_{p,q}\right)_{\curlyvee-\ell}}{[\curlyvee-\ell]_{p,q}}b_\curlyvee\leq 1.$$

Thus by using theorem 2.2, $F_u\left(f_\curlyvee(z)\right)\in\overline{\mathcal{H}_{p,q,\ell}^{\delta}}(p,q,\curlyvee,\ell,\delta)$ . ■

## 7. CONCLUSION

This article have been introduced a new class of multivalent harmonic functions defined by linear operator $\mathcal{H}_{p,q,\ell}^{\delta}(p,q,\curlyvee,\ell,\delta)$. Thus, some interesting outcome are studied and demonstrated such as the coefficient bound for the class $\mathcal{H}_{p,q,\ell}^{\delta}(p,q,\curlyvee,\ell,\delta)$. Consequently, some geometric properties have been investigated and analyzed for the consider functions, like the coefficient inequality, extreme points, convolution condition, convex linear combinations, and integral transforms in connection to the class $\mathcal{H}_{p,q,\ell}^{\delta}(p,q,\curlyvee,\ell,\delta)$. Finally, many problems remain opened for instance: extension of the obtained results to an operator applicable in partial deferential equations as in [16], [17].

## ACKNOWLEDGEMENT

The authors would like to thank the referees and editors for their constructive comments and suggestions.